\documentclass[12pt]{article}
\usepackage{latexsym,amsfonts,amssymb}
\usepackage[dvips]{color}
\usepackage{colordvi,multicol}

\def \cal{\mathcal}
\def\bb{\mathbb}

\newtheorem{thm}{Theorem}[section]
\newtheorem{cor}[thm]{Corollary}
\newtheorem{lem}[thm]{Lemma}

\newtheorem{defi}[thm]{Definition}
\newtheorem{rem}[thm]{Remark}

\begin{document}

\title{\bf Some inequalities and limit theorems under sublinear expectations}
\author{Ze-Chun Hu\thanks{Corresponding
author: Department of Mathematics, Nanjing University, Nanjing
210093, China\vskip 0cm E-mail address: huzc@nju.edu.cn}\,\, and Yan-Zhi Yang\\
 {\small Nanjing University}}
 \date{}
 \maketitle

\noindent{\bf Abstract}\quad In this note, we study inequality and limit theory under sublinear expectations.
We mainly prove Doob's inequality for submartingale and Kolmogrov's
inequality. By Kolmogrov's inequality, we obtain a special version of Kolmogrov's law of large numbers. Finally, we present a strong law of large numbers for independent and identically distributed random variables under one-order type moment condition.

\smallskip


\smallskip

%\noindent {\bf Mathematics Subject Classification (2000)}\quad Primary: 60J45;
%Secondary: 60G51

%\tableofcontents

\section{Introduction}
The classic strong law of large numbers  play an important role in the
development
of probability theory and its applications. One of essential ingredients
in these limit theorems is the additivity of the probabilities and the
expectations. However, many uncertain phenomena can not be well modeled
by using additive probabilities and additive expectations. Thus people
have used non-additive probabilities (called capacities) and nonlinear
expectations (for example Choquet integal/expectation, $g$-expectation) to
interpret and study these phenomena. Recently, motivated by the risk measures,
 superhedge pricing and modeling uncertain in finance, Peng \cite{Peng1}-
 \cite{Peng6}
 initiated the notion of independent and identically distributed (IID) random
 variables under sublinear expectations, and proved the weak law of large
 numbers and the central limit theorems. In \cite{Chen10}, Chen proved a
 strong law of large numbers for IID random variables
 under capacities induced by sublinear  expectations. In \cite{Hu12}, Hu
 presented three laws of large numbers for independent random
 variables without the requirement of identical distribution.

All the above existing results about the law of large numbers under
sublinear expectations need the moment condition of  $(1+\alpha)$-order
for some $\alpha>0$. But we know that in the classic additive
probability setting, if $X_1,X_2,\ldots$ is a sequence of IID random
variables with $E|X_1|<\infty$, then $(X_1+\cdots+X_n)/n$ converges
to $EX_1$ almost surely (a.s.). In virtue of this result, the
motivation of this note is to explore
 the law of large numbers for one sequence of IID random variables with one-order moment condition under
 sublinear expectations.

The rest of this note is organized as follows. In Section 2, we
recall some basic defintions, lemmas and theorems under sublinear
expectations. In Section 3,  We  prove Doob's inequality for
submartingale. In Section 4, we give Kolmogrov's inequality and some
applications to limit theory.

\section{Preliminary}
In this section, we present some basic definitions, lemmas and
theorems in the theory of sublinear expectations.

\begin{defi}\label{defi2.1}(see \cite{Peng1}-\cite{Peng6})
Let $(\Omega, \cal{F})$ be a measurable space and   $\cal{H}$ be a linear space of real valued functions defined on $\Omega$. We assume that $\cal{H}$ satisfies $c \in \cal{H}$ for any constant $c$ and $|X|\in \cal{H}$ if $X\in \cal{H}$.  $\cal{H}$ is considered as the space of our ``random variables". A nonlinear expectation $\bb{E}$ on $\cal{H}$ is a functional $\mathbb{E}: \cal{H} \mapsto \mathbb{R}$ satisfying the following properties: for all $X, Y\in \cal{H}$, we have \\
(a)  Monotonicity: if $X\geq Y$ then  $\bb{E}[X]\geq \bb{E}[Y]$.\\
(b)  Constant preserving: $\bb{E}[c]=c, \forall c\in \mathbb{R}.$\\
The triple $(\Omega, \cal{H}, \bb{E})$ is called a nonlinear expectation space. We are mainly concerned with sublinear expectation where the expectation $\bb{E}$ satisfies also\\
(c) Sub-additivity: $\bb{E}[X+Y]\leq  \bb{E}[X]+ \bb{E}[Y]$ .\\
(d)  Positive homogeneity: $ \bb{E}[\lambda X]=\lambda  \bb{E}[X] ,\forall \lambda\geq 0  .$ \\
If only (c) and (d) satisfied, $\bb{E}$ is called a sublinear functional.
\end{defi}

Given a sublinear expectations $\bb{E}$, let us denote the conjugate expectation $\hat{\bb{E}}$ of sublinear expectation $\bb{E}$ by
$$
\hat{\bb{E}}[X]:=-\bb{E}[-X],\quad \forall X\in \cal{H}.
$$
By $0=\bb{E}([X+(-X)]\leq \bb{E}[X]+\bb{E}[-X]$, we know that
for any $X\in \cal{H}$, $ \hat{\bb{E}}[X]\leq \bb{E}[X]$.

If for any $A\in\cal{F},I_A\in\cal{H}$, then we denote a pair
$(\bb{V}, v)$ of capacities by
$$
\bb{V}(A):=\bb{E}[I_{A}],\quad v(A):=\hat{\bb{E}}[I_{A}],\quad \forall A\in \cal{F}.
$$
It is easy to show that
$$
\bb{V}(A)+v(A^{c})=1,\quad \forall A\in \cal{F},
$$
where $A^{c}$ is the complement set of $A$.

\begin{defi}\label{defi2.2} (see \cite[Definition 2.2]{Chen10})
A set function $V: \cal{F} \rightarrow [0,1]$ is called a continuous capacity if it satisfies\\
(i) $V(\Phi)=0, V(\Omega)=1.$\\
(ii) $V(A)\leq V(B)$, whenever $ A\subset B$ and $ A, B\in\cal{F}$.\\
(iii) $V(A_{n})\uparrow V(A)$,  if $A_{n}\uparrow A $,  where $A_{n}, A\in \cal{F}$.\\
(iv) $V(A_{n})\downarrow V(A)$,  if $A_{n}\downarrow A $,  where
$A_{n}, A\in \cal{F}$.
\end{defi}

%
%\begin{ass}\label{ass2.5}
%Throughout this paper, we assume that $(\bb{V}, v)$ is pair of continuous capacities generated by sublinear expectation $\bb{E}$ and its conjugate expectation $\cal{E}$.
%\end{ass}

\begin{defi}\label{defi2.6} (see \cite[Definition 3]{Chen10})
\textbf{Independence:} Suppose that $Y_{1}, Y_2, \ldots, Y_n $ is a
sequence of random variables such that $Y_{i}\in \cal{H}$. Random
variable $Y_n $ is said to be independent of $X:=(Y_1, \ldots,
Y_{n-1})$ under $\bb{E}$, if for each measurable function $\varphi $
on ${\bb{R}}^n$ with $\varphi(X, Y_n)\in \cal{H}$ and $\varphi(x,
Y_n)\in \cal{H}$ for each $x\in {\bb{R}}^{n-1}$, we have
$$
\bb{E}[\varphi(X,Y_n)]=\bb{E}[\overline{\varphi}(X)],
$$
where $\overline{\varphi}(x):=\bb{E}[\varphi(x,Y_n)]$ and $\overline{\varphi}(X)\in \cal{H}$.

\textbf{Identical distribution:}  Random variables $X$ and $Y$ are
said to be identically distributed, denoted by  $X =^{\!\!\!\!d} Y$,
if for each $\varphi$ such that $\varphi(X), \varphi(Y)\in \cal{H}$,
$$
\bb{E}[\varphi(X)]=\bb{E}[\varphi(Y)] .
$$

\textbf{Sequence of IID random variables:} A sequence of random variables ${\{X_i\}}_{i=1}^{\infty}$ is said to be IID random variables, if $X_i=^{\!\!\!\!d} X_1$ and $X_{i+1}$ is independent of $Y:=(X_1, \cdots,X_i)$ for each $i\geq 1$.
\end{defi}

Since we will use some results about martingales in discrete time built  in \cite{CJP11}, now
we present some basic definitions and propositions in \cite{CJP11}.

Let $(\Omega,\cal{F})$ be a measurable space and $\{{\cal{F}}_t\}_{t\in \bb{N}}$ be a
discrete-time filtration on this space. Assume $\cal{F}={\cal{F}}_{\infty}={\cal{F}}_{\infty-}$
and ${\cal{F}}_0$ is trivial.  Let $m{\cal{F}}_t$ denote the space of ${\cal{F}}_t$-measurable
$\bb{R}\cup\{\pm \infty\}$-valued functions. The concepts of measurability, adaptedness, stopping times and the $\sigma$-algebra at a stopping time are identical to the classical case.

\begin{defi} (see \cite[Definition 2.1]{CJP11})\label{defi2.4}
Let $\cal{H}$ be a linear space of $\cal{F}$-measuable
$\bb{R}$-valued functions on $\Omega$ containing the constants. We
assume that $X\in \cal{H}$ implies that $|X|\in \cal{H}$ and
$I_{A}X\in \cal{H}$ for any $A\in \cal{F}$, and define
${\cal{H}}_t:=\cal{H}\cap m{\cal{F}}_t$.
\end{defi}

\begin{defi} (see \cite[Definition 2.2]{CJP11})\label{defi2.6}
A family of maps ${\cal{E}}_t :\cal{H}\rightarrow {\cal{H}}_t$ is called a ${\cal{F}}_t$-consistent nonlinear expectation, if for any $X, Y\in \cal{H}$, for all $s\leq t$,\\
(i) $X\geq Y $ implies ${\cal{E}}_t(X)\geq {\cal{E}}_t(Y)$;\\
(ii) ${\cal{E}}_s(Y)={\cal{E}}_s({\cal{E}}_t(Y))$;\\
(iii) ${\cal{E}}_t(I_{A}Y)=I_{A}{\cal{E}}_t(Y)$ for all $A \in {\cal{F}}_t$;\\
(iv) ${\cal{E}}_t(Y)=Y$ for all $Y\in{\cal{H}}_t$.\\
 A nonlinear expectation is called sublinear if it also satisfies\\
 (v) ${\cal{E}}_t(X+Y)\leq {\cal{E}}_t(X)+{\cal{E}}_t(Y)$;\\
 (vi) ${\cal{E}}_t(\lambda Y)={\lambda}^{+}{\cal{E}}_t(Y)+{\lambda}^{-}{\cal{E}}_t(-Y)$ for all $\lambda \in {\cal{H}}_t $ with $\lambda Y \in \cal{H}$.\\
 A nonlinear expectation is said to have the monotone continuity property (or Fatou property) if \\
 (vii) For any sequence $\{X_{i}\}$ in $\cal{H}$ such that $X_{i}(\omega)\downarrow 0$ for each $\omega$ , we have ${\cal{E}}_0(X_i) \rightarrow 0.$

  An ${\cal{F}}_t$-consistent sublinear expectation with the monotone continuity property will, for simplicity, be called an $\cal{SL}$-expectation. As ${\cal{F}}_{0}$ is trivial, one can equate ${\cal{E}}_0$ with a map $\cal{E}:\cal{H}\rightarrow \bb{R}$, satisfying the above properties.
  \end{defi}

\begin{lem} (Jensen's inequality)(see \cite[Lemma 2.3]{CJP11})\label{lem2.7} For any convex function $\varphi: \bb{R} \rightarrow \bb{R}$, any t, if $X$ and $\varphi(X)$ are both in $\cal{H}$, then
$$
\cal{E}_t(\varphi(X))\geq \varphi(\cal{E}_t(X)).
$$
\end{lem}

\begin{lem}(see \cite[Theorem 2.1]{Peng6})
An $\cal{SL}$-expectation has a representation
$$
\cal{E}(Y)=\sup\limits_{\theta\in \Theta}E_{\theta}[Y],
$$
where $\Theta$ is a collection of ($\sigma$-additive) probability measures on $\Omega$.
\end{lem}

\begin{defi}(see \cite[Definition 2.3]{CJP11})\label{q.s.}
We say that a statement holds quasi-surely (q.s.) if it holds except
on a set $N$ with $\cal{E}(I_N)=0$, or equivalently, if it holds
$\theta$-a.s. for all $\theta\in\Theta$.
\end{defi}

\begin{defi}(see \cite[Definition 2.4]{CJP11})\label{extention}
For any pair $(\cal{H},\cal{E})$, where $\cal{H}$ satisfies Definition \ref{defi2.6} and $\cal{E}$ is an $\cal{SL}$-expectation on $\cal{H}$, we can consistently extend our space to $({\cal{H}}^{ext},{\cal{E}}^{ext})$, where
\begin{eqnarray*}
&&{\cal{H}}^{ext}:=\{X\in m\cal{F}: \min\{E_{\theta}[X^{+}], E_{\theta}[X^{-}]\}<\infty\ \mbox{for all} \ \theta \in\Theta\},\\
&&{\cal{E}}^{ext}:=\sup\limits_{\theta\in\Theta}E_{\theta}[X].
\end{eqnarray*}
\end{defi}

\begin{defi}(see \cite[Definition 3.4]{CJP11})
A sequence $X_n\in \cal{H}^{ext}$ is said to  converge in capacity
to some $X_{\infty}\in \cal{H}^{ext}$ if, for any
$\varepsilon,\delta>0$, there exists an $N\in\bb{N}$ such that
$$
\cal{E}^{ext}(I_{\{|X_m-X_{\infty}|>\varepsilon\}})<\delta
$$
for all $m\geq N$.
\end{defi}

\begin{defi}(see \cite[Definition 2.5]{CJP11})
For a stopping time $T$,  define the expectation conditional on the $\sigma$-algebra ${\cal{F}}_T$ by
$$
{\cal{E}}_T(\cdot ;\omega):={\cal{E}}_{T(\omega)}(\cdot ;\omega),
$$
with the formal equivalence ${\cal{E}}_{\infty}(X)=X$.
\end{defi}

\begin{defi}(see \cite[Definition 3.1]{CJP11})
For $p \in [1, \infty)$, the map
$$
{|| \cdot ||}_p : X \mapsto (\cal{E}(|X|^p))^{1/p}
$$
forms a seminorm on $\cal{H}$. Similarly for $p=\infty $, where
$$
|| \cdot ||_{\infty} : X \mapsto \inf\limits_{x\in \bb{R}}\{x : \cal{E}(I_{\{|X|>x\}})=0\}.
$$

Define the space ${\cal{L}}^{p}(\cal{F})$ as the completion under ${||\cdot||}_{p}$ of the set
$$
\{X\in \cal{H}: {||X||}_{p} < \infty\}
$$
and then $L^{p}(\cal{F})$ as the equivalence classes of ${\cal{L}}^{p}$ modulo equality in ${||\cdot||}_{p}$.
\end{defi}

\begin{defi}(see \cite[Definition 3.2]{CJP11})
Consider $K\subset L^{1}$. $K$ is said to be uniformly integrable (u.i.)
if $\cal{E}(I_{\{|X|\geq c\}}|X|)$ converges to 0 uniformly in $X\in K$ as
$c\rightarrow \infty$.
\end{defi}

\begin{thm}(see \cite[Theorem 3.1]{CJP11})\label{thm2.13}
Suppose $K$ is a subset of $L^1$. Then $K$ is uniformly integrable
if and only if both\\
(i) $\{\cal{E}(|X|)\}_{X\in K}$ is bounded; and \\
(ii) For any $\epsilon>0$ there is a $\delta>0$ such that for all
$A\in\cal{F}$ with $\cal{E}(I_A)\leq \delta$ we have
$\cal{E}(I_A|X|)<\epsilon$ for all $X\in K$.
\end{thm}

By Theorem \ref{thm2.13}, we can easily get the following result.
\begin{cor}\label{cor2.14}
Suppose $K$ is a uniformly integrable family in $L^1$. Then its
closed convex hull in $L^1$ is also uniformly integrable.
\end{cor}

\begin{defi}(see \cite[Definition 3.3]{CJP11})
Let $L_{b}^{p}$ be the completion of the set of bounded functions $X \in \cal{H}$, under the norm ${|| \cdot ||}_p$ .
Note that $L_{b}^{p}\subset L^{p}$ .
\end{defi}

\begin{lem}(see \cite[Lemma 3.4]{CJP11})\label{lem2.13}
For each $p\geq 1$,
$$
L_{b}^{p}=\{X\in L^{p}: \lim_{n\rightarrow \infty}\cal{E}(|X|^{p}I_{\{|X|>n\}})=0\}.
$$
\end{lem}

\begin{thm}(see \cite[Theorem 3.2]{CJP11})\label{thm1}
Suppose $X_n$ is a sequence in $L_{b}^{1}$, and $X\in {\cal{H}}^{ext}$ .  Then the $X_n$ converge in $L^{1}$ norm to $X$ if and only if the collection ${\{X_n\}}_{n\in\bb{N}}$ is uniformly integrable and the $X_n$ converge in capacity to $X$.
Furthermore, in this case, the collection ${\{X_n\}}_{n\in\bb{N}}\cup \{X\}$ is also uniformly integrable and $X\in L_{b}^{1}$ .
\end{thm}

\begin{thm}
(The dominated  convergence theorem)
(see \cite[Theorem 3.11]{Chen12}) \label{thm2.16}
Suppose that $\{X_n\}$ is a sequence in $L^1$, and $|X_n|\leq Y, \forall n\geq 1$ with $Y\in L_{b}^{1}$ . If $X_n \rightarrow X$ q.s.  or $X_n$ converge to $X$ in capacity, then $X_n$ converge  to $X$ in $L^1$ norm.
\end{thm}
{\bf Proof.} Notice that $|X_n|\leq Y$ for all $n\geq 1$ and $Y\in
L_{b}^{1}$. Then by the definition of uniformly integrability and Lemma
\ref{lem2.13}, we have that $X_n\in L_{b}^{1}$ for all $n\geq 1$ and
$\{X_n\}$ is uniformly integrable. Then by Theorem \ref{thm1},
$X_n$ converge to $X$ in $L^1$ norm.\hfill\fbox

\begin{rem}
The condition $Y\in L_{b}^{1}$ in the above theorem can not be
weakened to be $Y\in L^1$. In fact, suppose that $Y\in L^1$. For $n\geq
1$, set $X_n:=|Y|I_{|Y|\geq n}$. Then $X_n$ converges to 0 q.s. But
by Lemma \ref{lem2.13},  $\cal{E}(|X_n-0|)\to 0$ if and only if
$Y\in L_b^1$.
\end{rem}

\begin{defi}(see \cite[Definition 4.1]{CJP11})
A process $X$ is called an $\cal{SL}$-martingale if it satisfies
$$
X_{s}={\cal{E}}_s(X_t)
$$
for all $s\leq t$, and $X_t \in L^{1}\cap m{\cal{F}}_t$  for all t. Similarly we define $\cal{SL}$-supermartingale and $\cal{SL}$-submartingale.
\end{defi}

%\begin{lem}
%Let $(X_t)$ be a $\cal{SL}$-supermartingale. Then $(-X_t)$ is a $\cal{SL}$-submartingale.\\
%(Note, the converse result does not hold.)
%\end{lem}

\begin{thm}(see \cite[Theorem 4.1]{CJP11})\label{thm2.20}
Let $X$ be an $\cal{SL}$-submartingale, and $S\leq T$ be bounded stopping times. Then $X_{S}\leq {\cal{E}}_{S}(X_T)$ .
\end{thm}

\section{Doob's inequality for submartingale}\setcounter{equation}{0}
In this section, suppose that $(\cal{E}_t)_{t\in \mathbb{N}}$ be a sublinear expectation defined by Definition \ref{defi2.6},
$\cal{E}$ denotes ${\cal{E}}_0$ and $(\bb{V},v)$ denotes the pair of  the capacities  generated by  $\cal{E}$ and its  conjugate expectation.

\begin{thm}\label{thm3.1}
Let ${\{X_{j},{\cal{F}}_j\}}_{j=1}^{n}$ be an $\cal{SL}$-submartingale, then for each $\lambda >0$,\\
(i) $\lambda \bb{V}\left(\max\limits_{1\leq j\leq n}X_{j}\geq \lambda\right)\leq \cal{E}\left(I_{\{\max\limits_{1\leq j\leq n}X_{j}\geq \lambda\}}X_n\right)\leq \cal{E}({X_n}^{+})\leq \cal{E}(|X_n|)$, where $X_n^+:=X_n\vee 0$;\\
(ii) $\lambda v\left(\min\limits_{1\leq j\leq n}X_{j}\leq -\lambda\right)\leq \cal{E}(X_n)-\cal{E}(X_1)+\cal{E}\left(-I_{\{\min\limits_{1\leq j\leq n}X_{j}\leq - \lambda\}}X_n\right)\leq \cal{E}(X_n)-\cal{E}(X_1)+\cal{E}({X_n}^{-})$, where $X_n^-:=-(X_n\wedge 0)$.
\end{thm}

To prove Theorem \ref{thm3.1}, we need the following lemma.

\begin{lem}(see \cite[Theorem 2.18(ii)(iv)]{Chen12})\label{lem3.2}\\
(i) Let $T$ be a finite stopping time. If
$A\in\cal{F}_T,X\in\cal{H}$, then $ \cal{E}_T(I_AX)=I_A\cal{E}_T(X).
$\\
(ii) For any two bounded stopping times $S,T$ with $S\leq T\leq S+1$ and $X\in \cal{H}$, we have $\cal{E}_{S}(\cal{E}_T(X))=\cal{E}_S(X)$.
\end{lem}

\noindent {\bf Proof of Theorem \ref{thm3.1}.} (i) Define $T:=\inf
\{j\geq 1, X_j\geq \lambda \}\wedge n$ and $A:=\{\max\limits_{1\leq
j\leq n}X_{j}\geq \lambda\}$. Then $T$ is a  stopping time with
$T\leq n$,  $A\in \cal{F}_T$ and  $X_{T}\geq \lambda$ on $A$.   By
Theorem \ref{thm2.20}, we have $ X_T\leq \cal{E}_T(X_n), $ and thus
$I_AX_T\leq I_A\cal{E}_T(X_n)$, which together with Lemma
\ref{lem3.2} implies that
\begin{eqnarray*}
\cal{E}(\lambda I_{A})&\leq& \cal{E}(I_{A}X_{T})\leq \cal{E}(I_A\cal{E}_T(X_n))\\
&=&\cal{E}(\cal{E}_T(I_AX_n))=\cal{E}(\cal{E}_{T\wedge 1}(\cal{E}_T(I_AX_n)))\\
&=&\cdots\\
&=&\cal{E}(\cal{E}_{T\wedge 1}(\cdots(\cal{E}_{T\wedge (n-1)}(\cal{E}_T(I_AX_n)))\cdots))\\
&=&\cal{E}(\cal{E}_{T\wedge 1}(\cdots(\cal{E}_{T\wedge (n-1)}(I_AX_n))\cdots))\\
&=&\cdots\\
&=&\cal{E}(\cal{E}_{T\wedge 1}(I_AX_n))\\
&=&\cal{E}(I_{A}X_{n})\leq \cal{E}({X_{n}}^{+}).
\end{eqnarray*}
Hence  $\lambda \bb{V}(A)\leq \cal{E}(I_{A}X_{n})\leq \cal{E}({X_{n}}^{+})\leq \cal{E}(|X_n|)$ .

(ii) Let $T:=\inf\{1\leq j\leq n,X_{j}\leq -\lambda\}\wedge (n+1), M_{k}:=\{\min\limits_{1\leq j\leq k}X_{j}\leq -\lambda\}=\{T\leq k\}$. Then $T$ is a  stopping time satisfying $1\leq T\leq n+1$. By Theorem \ref{thm2.20}, we have
$X_1\leq \cal{E}_1(X_{T\wedge n})$. By Definition \ref{defi2.6} and the properties of sublinear expectation, we have
\begin{eqnarray*}
\cal{E}(X_{1})&\leq & \cal{E}(\cal{E}_1(X_{T\wedge n}))=\cal{E}(X_{T\wedge n})\\
&=&\cal{E}(I_{\{T\leq n\}}X_{T}+I_{{M_{n}}^{c}}X_{n})\\
&\leq &\cal{E}(-\lambda I_{M_{n}}) +\cal{E}(I_{{M_{n}}^{c}}X_{n})\\
&=&-\lambda v(M_{n})+\cal{E}(X_{n}-I_{M_{n}}X_{n})\\
&\leq & -\lambda v (M_{n})+\cal{E}(X_{n})+\cal{E}(-I_{M_{n}}X_{n}).
\end{eqnarray*}
It follows that the two inequalities in (ii) hold.\hfill\fbox

\begin{cor}\label{c1}
Let ${\{X_{j},{\cal{F}}_j\}}_{j=1}^{n}$ be an $\cal{SL}$-martingale. Then for each $\lambda >0$, we have\\
(i) $\bb{V}\left(\max\limits_{1\leq j\leq n}|X_{j}|\geq \lambda\right)\leq \frac{1}{\lambda}\cal{E}\left(I_{\{\max\limits_{1\leq j\leq n}|X_{j}|\geq \lambda\}}|X_{n}|\right)\leq\frac{1}{\lambda}\cal{E}(|X_{n}|).
$\\
(ii) If $\cal{E}(X_j^2)<\infty, \forall 1 \leq j\leq n$, then
$
\bb{V}\left(\max\limits_{1\leq j\leq n}|X_{j}|\geq \lambda\right)\leq \frac{1}{\lambda^{2}}\cal{E}(X_{n}^{2}).
$
\end{cor}
{\bf Proof.} By Jensen's inequality (see Lemma \ref{lem2.7}), we
know that ${\{|X_{j}|,{\cal{F}}_j\}}_{j=1}^{n}$ in (i)
(respectively, ${\{X_{j}^2,{\cal{F}}_j\}}_{j=1}^{n}$ in (ii)) is an
$\cal{SL}$-submartingale. Then the results follows from Theorem
\ref{thm3.1}(i).\hfill\fbox

\section{Kolmogrov's inequality and its applications}\setcounter{equation}{0}

Let $(\Omega,\cal{F})$ be a measurable space and $\{{\cal{F}}_t\}_{t\in \bb{N}}$ be a
discrete-time filtration on this space. Assume $\cal{F}={\cal{F}}_{\infty}={\cal{F}}_{\infty-}$
and ${\cal{F}}_0$ is trivial.  Let $m{\cal{F}}_t$ denote the space of ${\cal{F}}_t$-measurable
$\bb{R}\cup\{\pm \infty\}$-valued functions. Let $\cal{H},\cal{H}_t$ be defined by Definition \ref{defi2.4} and $\{\cal{E}_t\}_{t\in\bb{N}}$ be an $\cal{SL}$-expectation defined by Definition \ref{defi2.6}. Assume that the capacity $\bb{V}$ associated with $\cal{E}_0$ (i.e. $\cal{E}$) is a continuous capacity defined by Definition \ref{defi2.2}.

\subsection{Some preparations}
\begin{lem}\label{lem4.1}
Suppose that $X,Y\in \cal{H}_t$ for some $t\geq 0$. If  $\cal{E}((X+c)I_A)=\cal{E}((Y+c)I_A)$ for any $A\in\cal{F}_t$ and $X+c\geq 0,Y+c\geq 0$ q.s. for some constant $c$, then $X=Y$ q.s.
\end{lem}
{\bf Proof.} For simplicity, we assume that $c=0$. Suppose that $X=Y$ q.s. doesn't hold. Then $\bb{V}(\{X>Y\}\cup\{X<Y\})>0$. Without loss of generality, assume that $\bb{V}(\{X>Y\})>0$. Then there exists two constants $r,s$ with $r>s\geq 0$ such that $\bb{V}(\{X\geq r>s\geq Y\})>0$. Let $A:=\{X\geq r>s\geq Y\}$. Then $A\in\cal{F}_t$ and
$$
\cal{E}(XI_A)\geq r\bb{V}(A)>s\bb{V}(A)\geq \cal{E}(YI_A),
$$
which contradicts the assumption. Hence  $X=Y$ q.s.\hfill\fbox

\bigskip

\begin{defi}\label{defi4.2}
An element $X\in\cal{H}$ is said to be independent of $\cal{F}_n$ for some $n\in\bb{N}$ if $X$ is independent of $I_A$ under $\cal{E}$ for any $A\in\cal{F}_n$.
\end{defi}

\begin{lem}\label{lem4.3}
Let  $X\in\cal{H}$ be independent of $\cal{F}_n$ for some $n\in\bb{N}$. Then $\cal{E}_n(X)=\cal{E}(X)$ q.s. if one of the following two conditions holds:\\
(i) $X$ is lower bounded;\\
(ii) $X\in L_b^1$.
\end{lem}
{\bf Proof.} (i) By the definition of independence, we know that for any constant $c$, any $A\in\cal{F}_n$, $X+c$ is independent of $I_A$. Further, for any constant $c$, we have that $\cal{E}(X+c)=\cal{E}(X)+c$ and $\cal{E}_n(X+c)=\cal{E}_n(X)+c$. Hence without loss of generality, we assume that $X\geq 0$. Then $\cal{E}(X)\geq 0,\cal{E}_n(X)\geq 0$, and both are in $\cal{H}_n$. By the independence,  for any $A\in\cal{F}_n$,
\begin{eqnarray}\label{lem4.3-a}
\cal{E}(I_AX)=\cal{E}(I_A\cal{E}(X)).
\end{eqnarray}
By the properties of $\cal{E}_n$, for any $A\in\cal{F}_n$,
\begin{eqnarray}\label{lem4.3-b}
\cal{E}(I_A\cal{E}_n(X))=\cal{E}(\cal{E}_n(I_AX))=\cal{E}(I_AX).
\end{eqnarray}
By (\ref{lem4.3-a}) and (\ref{lem4.3-b}), we have that for any $A\in\cal{F}_n$,
$$
\cal{E}(I_A\cal{E}_n(X))=\cal{E}(I_A\cal{E}(X)),
$$
which together with Lemma \ref{lem4.1} implies that $\cal{E}_n(X)=\cal{E}(X)$ q.s.

(ii) For any $m\in\bb{N}$, define
$$
X_m:=XI_{|X|\leq m},\ Y_m:=\cal{E}_n(X_m).
$$
Since $X\in L_b^1$, by the dominated convergence theorem (see Theorem \ref{thm2.16}), we have
\begin{eqnarray}\label{lem4.3-c}
\lim_{m\to\infty}\cal{E}(|X_m-X|)=0.
\end{eqnarray}
It follows that
\begin{eqnarray}\label{lem4.3-d}
\lim_{m\to\infty}\cal{E}(X_m)=\cal{E}(X).
\end{eqnarray}
By (\ref{lem4.3-c}) and \cite[Lemma 3.1]{CJP11}, we have
\begin{eqnarray}\label{lem4.3-e}
\lim_{m\to\infty}\cal{E}(|\cal{E}_n(X_m)-\cal{E}_n(X)|)=0.
\end{eqnarray}
Since for any $m\in\bb{N}$, $X_m$ is lower bounded, then by (i) we have
\begin{eqnarray}\label{lem4.3-f}
\cal{E}_n(X_m)=\cal{E}(X_m)\ q.s.
\end{eqnarray}
By (\ref{lem4.3-d})-(\ref{lem4.3-f}) and the triangle inequality, we obtain that $\cal{E}(|\cal{E}_n(X)-\cal{E}(X)|)=0$, and so $\cal{E}_n(X)=\cal{E}(X)$ q.s.\hfill\fbox

\subsection{Kolmogorov's inequality}

\begin{thm}\label{thm4.1}
 Let $X_1, X_2,\ldots$ be a sequence of  random variables on the sublinear expectation space
$(\Omega, \cal{H},\cal{E})$ such that for any $i=1,2,\ldots,X_i\in\cal{H}_i$, $X_i$ is independent of $\cal{F}_{i-1}$,  $X_i\in L_b^1$, and $\cal{E}(-X_i)=-\cal{E}(X_i)$. For  any $j\in \bb{N}$, define
$S_j=\sum_{i=1}^j(X_i- \cal{E}(X_i))$. Then for any $n\in \bb{N}$
and $\varepsilon >0$, we have
$$
\bb{V}\left(\max\limits_{1\leq j\leq n}|S_{j}|\geq
\varepsilon\right)\leq \frac{1}{\varepsilon^{2}}\cal{E}(S_n^2)= \frac{1}{\varepsilon^{2}}
\sum\limits_{i=1}^{n}\cal{E}((X_{i}-\cal{E}(X_i))^{2}).
$$
\end{thm}
{\bf Proof.} Without loss of generality, we can assume that
$\cal{E}(X_j)=\cal{E}(-X_j)=0,\forall j\in \bb{N}$.   For
$0\leq n\leq m$, we consider $ \cal{E}_n(S_m)$. If $n=m$,
then $\cal{E}_n(S_m)= \cal{E}_n(S_n)=S_n$. If
$n<m$, then by Definition \ref{defi2.6} and Lemma \ref{lem4.3}, we have
\begin{eqnarray*}
\cal{E}_n(S_m)&=&
\cal{E}_n(X_1+\cdots+X_n+X_{n+1}+\cdots+X_m)\\
&=&X_1+\cdots+X_{n}+\cal{E}_n(X_{n+1}+\cdots+X_m)\\
&=&X_1+\cdots+X_{n}+\cal{E}_n(\cal{E}_{m-1}(X_{n+1}+\cdots+X_m))\\
&=&X_1+\cdots+X_{n}+\cal{E}_n(X_{n+1}+\cdots+X_{m-1}+\cal{E}_{m-1}(X_m))\\
&=&X_1+\cdots+X_{n}+\cal{E}_n(X_{n+1}+\cdots+X_{m-1}+\cal{E}(X_m))\\
&=&X_1+\cdots+X_{n}+\cal{E}_n(X_{n+1}+\cdots+X_{m-1})\\
&=&\cdots\\
&=&X_1+\cdots+X_{n}+\cal{E}_n(X_{n+1})\\
&=&X_1+\cdots+X_{n}+\cal{E}(X_{n+1})\\
&=&X_1+\cdots+X_{n}=S_n.
\end{eqnarray*}
Hence $\{S_j\}_{j\in \bb{N}}$ is an $\cal{SL}$-martingale.
By Corollary \ref{c1}, we have
\begin{eqnarray*}
&&\bb{V}\left(\max\limits_{1\leq j\leq n}|S_{j}|\geq \epsilon\right)
\leq \frac{1}{\varepsilon^{2}}\cal{E}(|S_n|^{2}).
\end{eqnarray*}
For any $1\leq i<j$, by Definition \ref{defi2.6} and Lemma \ref{lem4.3},  we have
\begin{eqnarray}\label{thm4.1-a}
\cal{E}(X_{i}X_{j})&=&\cal{E}(\cal{E}_i(X_{i}X_{j}))=\cal{E}(X_{i}^+\cal{E}_i(X_{j})+X_{i}^-\cal{E}_i(-X_{j}))\nonumber\\
&=&\cal{E}(X_{i}^+\cal{E}(X_{j})+X_{i}^-\cal{E}(-X_{j}))=0.
\end{eqnarray}
Similarly, we have
\begin{eqnarray}\label{thm4.1-b}
\cal{E}(-X_{i}X_{j})&=&\cal{E}((-X_{i})^+\cal{E}(X_{j})+(-X_{i})^-\cal{E}(-X_{j}))=0.
\end{eqnarray}
Then
by (\ref{thm4.1-a}), (\ref{thm4.1-b}), the properties of sublinear expectations,  and the independence of
$\{X_j\}_{j\in\bb{N}}$, we have
\begin{eqnarray*}
\cal{E}(|S_n|^{2})&=&\cal{E}((X_{1}+\cdots+X_n)^{2})\\
&=&\cal{E}\left(\sum_{i=1}^nX_i^2+ 2\sum_{1\leq i<j\leq n}
X_iX_j\right)\\
&=& \sum\limits_{i=1}^{n}\cal{E}(X_i^{2}).
\end{eqnarray*}
\hfill\fbox

\subsection{Some applications}

\begin{thm}\label{thm4.2}
Let $X_1, X_2,\ldots$ be a sequence of  random variables on the sublinear expectation space
$(\Omega, \cal{H},\cal{E})$ such that for any $i=1,2,\ldots,X_i\in\cal{H}_i$, $X_i$ is independent of $\cal{F}_{i-1}$, $X_i\in L_b^1$, and $\cal{E}(-X_i)=-\cal{E}(X_i)$. If
$\sum\limits_{i=1}^{\infty}\cal{E}((X_{i}-\cal{E}(X_i))^2)<\infty$,
then $\sum\limits_{i=1}^{\infty}(X_i-\cal{E}(X_i))$ is q.s.
convergent.
\end{thm}
{\bf Proof.} Without loss of generality, we  assume that
$\cal{E}(X_i)=\cal{E}(-X_i)=0,\forall i\in \bb{N}$. For any $n\in \mathbb{N}$, let $S_n:=X_1+\cdots+X_n.$ Denote
$A:=\{\omega\in\Omega: S_{j}(\omega)-S_{k}(\omega)\nrightarrow 0 \ \ \mbox{as}\ j,k\rightarrow \infty\}$.
By the definition of q.s. convergence, we know that $S_n$ is q.s. convergent if and only if $\bb{V}(A)=0$. The set $A$ can be expressed by
$
A=\cup_{k=1}^{\infty}\cap_{n=1}^{\infty}\cup_{j,k=n}^{\infty}\{|S_{j}-S_{k}|\geq \frac{1}{k}\}.
$
If $\bb{V}(A)=0$, then by the monotone property of $\bb{V}$, we have that for any $k\in \mathbb{N}$,
\begin{eqnarray}\label{thm4.2-a}
\bb{V}\left(\cap_{n=1}^{\infty}\cup_{j,k=n}^{\infty}\{|S_{j}-S_{k}|\geq \frac{1}{k}\}\right)=0.
\end{eqnarray}
Conversely, if for any $k\in \mathbb{N}$, (\ref{thm4.2-a}) holds, then by the continuity and sub-additivity of $\bb{V}$, we have
\begin{eqnarray}\label{thm4.2-b}
\bb{V}(A)&=&\lim_{m\to\infty}\bb{V}\left(\cup_{k=1}^m\cap_{n=1}^{\infty}\cup_{j,k=n}^{\infty}\{|S_{j}-S_{k}|\geq \frac{1}{k}\}\right)\nonumber\\
&\leq &\limsup_{m\to\infty}\sum_{k=1}^m\bb{V}\left(\cap_{n=1}^{\infty}\cup_{j,k=n}^{\infty}\{|S_{j}-S_{k}|\geq \frac{1}{k}\}\right)=0.
\end{eqnarray}
Hence $\bb{V}(A)=0$ if and only if for any $k\in \mathbb{N}$, (\ref{thm4.2-a}) holds, or equivalently, if and only if for any $\varepsilon>0$,
\begin{eqnarray}\label{thm4.2-c}
\bb{V}\left(\cap_{n=1}^{\infty}\cup_{j,k=n}^{\infty}\{|S_{j}-S_{k}|\geq \varepsilon\}\right)=0.
\end{eqnarray}
By the continuity of $\bb{V}$, (\ref{thm4.2-c}) holds if and only if
\begin{eqnarray}\label{thm4.2-d}
\lim_{n\to\infty}\bb{V}\left(\cup_{j,k=n}^{\infty}\{|S_{j}-S_{k}|\geq \varepsilon\}\right)=0.
\end{eqnarray}
If $j>n,k>n$, then we have
\begin{eqnarray}\label{thm4.2-e}
\{|S_{j}-S_{k}|\geq \varepsilon\}\subset \{|S_{j}-S_{n}|\geq \frac{\varepsilon}{2}\}\cup\{|S_{k}-S_{n}|\geq \frac{\varepsilon}{2}\}.
\end{eqnarray}
By (\ref{thm4.2-e}), we know that (\ref{thm4.2-d}) holds for any $\varepsilon>0$ if and only if
\begin{eqnarray}\label{thm4.2-f}
\lim_{m\to\infty}\bb{V}(\cup_{k=1}^{\infty}\{|S_{m+k}-S_{k}|\geq \varepsilon\})=0
\end{eqnarray}
for any $\varepsilon>0$.

By the continuity of  the capacity $\bb{V}$, Theorem \ref{thm4.1},
and the condition
that $\sum\limits_{i=1}^{\infty}\cal{E}(X_{i}^2)<\infty$,
we have
\begin{eqnarray*}
&&\bb{V}(\cup_{k=1}^{\infty}\{|S_{m+k}-S_{k}|\geq \varepsilon\})\\
&&= \lim\limits_{n\rightarrow \infty}\bb{V}
(\cup_{k=1}^{n}\{|S_{m+k}-S_{k}|\geq \varepsilon\})\\
&&=\lim\limits_{n\rightarrow{\infty}}\bb{V}(\max\limits_{1\leq k\leq n}|S_{m+k}-S_{m}|\geq \varepsilon)\\
&&\leq \lim_{n\rightarrow{\infty}}\frac{1}{\varepsilon^{2}}
\sum\limits_{k=1}^{n}\cal{E}(X_{m+k}^2)\\
&&=\frac{1}{\varepsilon^{2}}
\sum_{k=m+1}^{\infty}\cal{E}(X_k^2)\rightarrow 0\ \ \mbox{as}\ m
\rightarrow \infty.
\end{eqnarray*}
The proof is complete. \hfill\fbox

\begin{thm}\label{thm4.3}
(\textbf{A special version of Kolmogorov's strong law of large
numbers}) Let $X_1, X_2,\ldots$ be a sequence of  random variables on the sublinear expectation space
$(\Omega, \cal{H},\cal{E})$ such that for any $i=1,2,\ldots,X_i\in\cal{H}_i$, $X_i$ is independent of $\cal{F}_{i-1}$, $X_i\in L_b^1$, and $\cal{E}(-X_i)=-\cal{E}(X_i)$. Let  $\{b_n\}$ be a sequence of increasing
positive numbers with $b_n \rightarrow \infty $ as $n \rightarrow
\infty$. If
$\sum\limits_{n=1}^{\infty}\frac{(\cal{E}(X_{n}-\cal{E}(X_n))^2)}{b_{n}^{2}}
<\infty$,  then $\frac{\sum_{i=1}^n(X_i-\cal{E}(X_i))}{b_{n}}
\rightarrow 0$ q.s. as $n\to\infty$.
\end{thm}
{\bf Proof.} By the assumption, we have
$\sum\limits_{n=1}^{\infty}\cal{E}
((\frac{X_n-\cal{E}(X_n)}{b_n})^2)<\infty$. Then by Theorem
\ref{thm4.2}, we know that
$\sum\limits_{i=1}^{\infty}(\frac{X_{i}}{b_i}-\cal{E}(\frac{X_{i}}{b_i}))$
 i.e. $\sum\limits_{i=1}^{\infty}\frac{X_{i}-\cal{E}(X_i)}{b_i}$
 is q.s. convergent. By  Kronecker's lemma, we obtain
$$
\frac{\sum_{i=1}^n(X_i-\cal{E}(X_i))}{b_{n}}=
\frac{1}{b_n}\sum_{i=1}^nb_i\frac{X_i-\cal{E}(X_i)}
{b_{i}}\to 0\ \ \mbox{q.s. as}\ \ n\to\infty.
$$
\hfill\fbox

\bigskip

To state the final result of this note, we give an assumption and Kolmogorov's 0-1 law (see \cite{Chen12}).

 {\bf Assumption (A):} For any sequence $\{X_n,n\in\bb{N}\}$ of random variables and
$B\in \sigma(X_n,n\in\bb{N})$, there is a sequence $\{B_n,n\in\bb{N}\}$ such that $B_n\in\sigma(X_1,\ldots,X_n)$ and $\lim\limits_{n\to\infty}\bb{V}(B_n\triangle B)=0$, where $B_n\triangle B:=(B_n\backslash B)\cup(B\backslash B_n)$.

\begin{rem}
Suppose that the sublinear expectation $\cal{E}$ can be expressed by
$$
\cal{E}(f)=\sup_{P\in\cal{P}}E_{P}(f),
$$
where $\cal{P}$ is a finite set of probability measures. Denote $\cal{P}:=\{P_1,\ldots,P_m\}$ and define
$$
\mu(A):=\sum_{k=1}^mP_k(A),\forall A\in\cal{F}.
$$
Then $\mu$ is a finite measure on $(\Omega,\cal{F})$. By the classic result in measure theory, for any $B\in\sigma(X_n,n\in\bb{N})$,
there is a sequence $\{B_n,n\in\bb{N}\}$ such that $B_n\in\sigma(X_1,\ldots,X_n)$ and $\lim\limits_{n\to\infty}\mu(B_n\triangle B)=0$. It follows that
\begin{eqnarray*}
\bb{V}(B_n\triangle B)=\cal{E}(I_{B_n\triangle B})=\sup_{P\in\cal{P}}P(B_n\triangle B)\leq \mu(B_n\triangle B)\to 0.
\end{eqnarray*}
\end{rem}

\begin{lem} (\textbf{Kolmogorov's 0-1 law})(see \cite[Theorem 4.15]{Chen12})\label{0-1 law}
Suppose that Assumption (A) holds and $\{X_n,n\in\bb{N}\}$ is a sequence of independent random variables. Then for any set $B$ in the tail $\sigma$-algebra $\cap_{n=1}^{\infty}\sigma(X_k,k\geq n)$, we have $\bb{V}(B)=0$ or 1.
\end{lem}

\begin{thm}\label{main}
Let $X_1, X_2,\ldots$ be a sequence of  random variables on the sublinear expectation space
$(\Omega, \cal{H},\cal{E})$ such that for any $i=1,2,\ldots,X_i\in\cal{H}_i$, $X_i$ is independent of $\cal{F}_{i-1}$, and they have the common distribution with
 $\sum\limits_{m=1}^{\infty}
 \cal{E}\left(|X_1|I_{\{m-1<|X_1|\leq m\}}\right)
 <\infty$. Then\\
 (i) if $\cal{E}(X_1)=-\cal{E}(-X_1)$, and  $\cal{E}(-X_1I_{\{|X_1|\leq n\}})=-\cal{E}(X_1I_{\{|X_1|\leq n\}}), \forall n\in\bb{N}$,
 then $S_n/n$ converges to $\cal{E}(X_1)$  q.s.;\\
 (ii) if $S_n/n$ converges to some $X\in\cal{H}$  q.s. and Assumption (A) holds, then $\cal{E}(X_1)=-\cal{E}(-X_1)$ and $X=\cal{E}(X_1)$ q.s.
    \end{thm}

 \begin{rem}
 (i) If $\{X_i,i\in\bb{N}\}$ is a sequence of IID random variables such that $\cal{E}(|X_1|^{1+\alpha})<\infty$ for some $\alpha>0$, then by \cite[Theorem 1]{Chen10},  $S_n/n$ converges to some $X\in\cal{H}$  q.s. if and only if $\cal{E}(X_1)=-\cal{E}(-X_1)$. In this case $X=\cal{E}(X_1)$ q.s.\\
 (ii) If $\cal{E}$ is a linear expectation, then the condition
$\sum_{m=1}^{\infty}\cal{E}(|X_1|I_{[m-1<|X_1|\leq m]})<\infty$ is
just $\cal{E}(|X_1|)<\infty$, and it holds that  $\cal{E}(-X_1I_{\{|X_1|\leq n\}})=-\cal{E}(X_1I_{\{|X_1|\leq n\}}), \forall n\in\bb{N}$.\\
(iii) As to the condition that $\cal{E}(-X_1I_{\{|X_1|\leq n\}})=-\cal{E}(X_1I_{\{|X_1|\leq n\}}), \forall n\in\bb{N}$, we give some sufficient conditions as follows:

(a) If $X_1$ is maximal distributed and $-\cal{E}(-X_1)=\cal{E}(X_1)$, then we can get that $\cal{E}(-X_1I_{\{|X_1|\leq n\}})=-\cal{E}(X_1I_{\{|X_1|\leq n\}}), \forall n\in\bb{N}$.

(b) Let the sublinear expectation $\cal{E}$  be expressed by $\cal{E}(\cdot)=\sup_{p\in \cal{P}}E_p(\cdot)$, where $\cal{P}$ is a family of probability measures and $E_p$ stands for the corresponding
(linear) expectation.  Suppose that $X_1$ is a symmetric random variable with respect to each probability measure $p\in\cal{P}$. Then for any $p\in\cal{P}$ and any $n\in \bb{N}$, we have
$E_p(X_1I_{\{|X_1|\leq n\}})=E_p(-X_1I_{\{|X_1|\leq n\}})=0$, and thus $\cal{E}(X_1I_{\{|X_1|\leq n\}})=\cal{E}(-X_1I_{\{|X_1|\leq n\}})=0$. \\
(iv) A natural question is

{\bf Under the condition that $\sum\limits_{m=1}^{\infty}
 \cal{E}\left(|X_1|I_{\{m-1<|X_1|\leq m\}}\right)
 <\infty$ or weaker condition (for example, $X_1\in L_b^1$), can we obtain the same strong law of large numbers
 with the one in \cite[Theorem 1]{Chen10} ?}

\noindent Unfortunately, by our method, we can not obtain that result.
It needs some new idea.

 \end{rem}

\noindent{\bf Proof of Theorem \ref{main}.}
 At first,
we claim that $X_1\in L_b^1$. In fact,
\begin{eqnarray*}
\cal{E}\left(|X_1|I_{\{|X_1|>n\}}\right)&=&
\cal{E}\left(|X_1|\left(\sum_{m=n}^{\infty}I_{\{m<|X_1|\leq m+1\}}\right)\right)\\
&\leq&\sum_{m=n}^{\infty}\cal{E}\left(|X_1|I_{\{m<|X_1|\leq m+1\}}\right)\\
&\to& 0\ \ \mbox{as}\ \  n\to\infty.
\end{eqnarray*}
Hence by Lemma \ref{lem2.13}, $X_1\in L_b^1$.

(i) Without loss of generality, we assume that
$\cal{E}(X_1)=0$.

{\it Step 1.} For $n\in\bb{N}$, let $Y_n:=X_nI_{\{|X_n|\leq n\}}$.
We claim that it's enough to show that $\frac{\sum_{i=1}^{n}Y_i}{n}
\rightarrow 0$ q.s.

By the properties of sublinear expectation
$\cal{E}$ and the assumption, we have
\begin{eqnarray*}
&&\sum_{k=1}^{\infty}\bb{V}(|X_1|>k)
=\sum_{k=1}^{\infty}\cal{E}\left(
\sum_{m=k}^{\infty}I_{\{m<|X_1|\leq m+1\}}\right)\\
&&\leq \sum_{k=1}^{\infty}\sum_{m=k}^{\infty}
\cal{E}(I_{\{m<|X_1|\leq m+1\}})\\
&&\leq \sum_{k=1}^{\infty}\sum_{m=k}^{\infty}
\cal{E}\left(\frac{|X_1|}{m}I_{\{m<|X_1|\leq m+1\}}\right)\\
&&=\sum_{m=1}^{\infty}\sum_{k=1}^{m}\frac{1}{m}
\cal{E}(|X_1|I_{\{m<|X_1|\leq m+1\}})\\
&&=\sum_{m=1}^{\infty}\cal{E}(|X_1|I_{\{m\leq|X_1|<m+1\}})<\infty.
\end{eqnarray*}
It follows that
$$
\sum_{n=1}^{\infty}\bb{V}(X_n\neq Y_n)=
\sum_{n=1}^{\infty}\bb{V}(|X_n|> n)=
\sum_{n=1}^{\infty}\bb{V}(|X_1|> n)<\infty.
$$
By Borel-Cantelli Lemma (see \cite[Lemma 2]{Chen10}), we have
$$
\bb{V}(\{X_n\neq Y_n\}, \mbox{i.o.})=0,
$$
where i.o. means ``infinitely often". Hence
$$
\frac{\sum_{i=1}^nX_i}{n}\rightarrow 0\ q.s.
\Longleftrightarrow \frac{\sum_{i=1}^nY_i}{n}\rightarrow 0\ q.s.
$$

{\it Step 2.} By the fact that $X_1\in L_b^1$ and the dominated
convergence theorem (see Theorem \ref{thm2.16}), we have
$$\cal{E}(Y_n)=\cal{E}(X_nI_{\{|X_n|\leq n\}})=\cal{E}(X_1I_{\{|X_1|\leq n\}})\rightarrow \cal{E}(X_1)=0\ \ \mbox{as}\ n\to \infty.
$$
Hence in order to prove that $\frac{\sum_{i=1}^nY_i}{n}\rightarrow
0$ q.s., it's enough to show that
$\frac{1}{n}\sum_{i=1}^{n}(Y_i-\cal{E}(Y_i))\rightarrow 0$ q.s.
Further, by Kronecker's lemma, it's enough to show that
$\sum_{i=1}^{\infty}\frac{Y_i-\cal{E}(Y_i)}{i}$ is q.s. convergent.
By Theorem \ref{thm4.2}, we need only to prove that
\begin{eqnarray}\label{thm4.4-a}
\sum_{n=1}^{\infty}\cal{E}\left(\left(\frac{Y_n-\cal{E}(Y_n)}{n}\right)^2\right)<\infty.
\end{eqnarray}

{\it Step 3.} Prove (\ref{thm4.4-a}). By Cauchy-Schwartz inequality,
the properties of sublinear expectations and the assumption, we have
\begin{eqnarray*}
\sum_{n=1}^{\infty}\cal{E}\left(\left(\frac{Y_n-\cal{E}(Y_n)}{n}\right)^2\right)
&\leq &\sum_{n=1}^{\infty}\frac{4}{n^2}\cal{E}(Y_n^2)
=\sum_{n=1}^{\infty}\frac{4}{n^2}\cal{E}({X_n}^2I_{\{|X_n|\leq n\}})\\
&\leq &\sum_{n=1}^{\infty}\frac{4}{n^2}\sum_{m=1}^{n}\cal{E}(X_1^2I_{\{m-1<|X_1|\leq m\}})\\
&=&\sum_{m=1}^{\infty}\cal{E}({X_1}^2I_{\{m-1<|X_1|\leq m\}})\sum_{n=m}^{\infty}\frac{4}{n^2}\\
%&\leq &8 \sum\limits_{m=1}^{\infty}\mathbb{E}[{X_1}^2I_{[m-1<|X_1|\leq m]}]\sum\limits_{n=m}^{\infty}\frac{1}{n(n+1)}\\
&\leq&\sum\limits_{m=1}^{\infty}\frac{4}{m}\cal{E}({X_1}^2I_{\{m-1<|X_1|\leq m\}})\\
&\leq &4 \sum_{m=1}^{\infty}\cal{E}(|X_1|I_{\{m-1<|X_1|\leq
m\}})<\infty.
\end{eqnarray*}

(ii) Suppose that $S_n/n$ converges to some $X\in
\cal{H}$ q.s. Since $\{X_i,i\in\bb{N}\}$ is a sequence of IID, by
the fact that $X_1\in L_b^1$, we know that $\{X_i,i\in \bb{N}\}$ is
uniformly integrable. Then by Corollary \ref{cor2.14}, we know that
$\{S_n/n,n\in\bb{N}\}$ is uniformly integrable. By Theorem
\ref{thm1}, we know that $S_n/n$ converges to $X$ in $L^1$-norm,
i.e.
\begin{eqnarray}\label{4.2}
\lim_{n\to\infty}\cal{E}(|\frac{S_n}{n}-X|)=0.
\end{eqnarray}
It follows that
\begin{eqnarray}\label{4.3}
\cal{E}(X)=\lim_{n\to\infty}\cal{E}(\frac{S_n}{n})=\cal{E}(X_1).
\end{eqnarray}
Similarly, by (\ref{4.2}), we have
\begin{eqnarray}\label{4.4}
-\cal{E}(-X)=-\lim_{n\to\infty}\cal{E}(-\frac{S_n}{n})=-\cal{E}(-X_1).
\end{eqnarray}

By Lemma \ref{0-1 law}, we know that  $X=c$ q.s. for some constant $c$. Thus
\begin{eqnarray}\label{4.5}
-\cal{E}(-X)=c=\cal{E}(X).
\end{eqnarray}
By (\ref{4.3})-(\ref{4.5}), we obtain that $\cal{E}(X_1)=-\cal{E}(-X_1)$ and $X=\cal{E}(X_1)$ q.s.
\hfill\fbox

\bigskip

 \noindent {\bf\large Acknowledgments} \vskip 0.1cm  \noindent

We acknowledge the helpful suggestions and comments of an anonymous referee,
which improved the presentation of this paper. We are grateful for the support of
the NNSFC, Jiangsu Province Basic Research Program (Natural Science
Foundation) (Grant No. BK2012720).

%\hypersetup{urlcolor=blue}

%\renewcommand{\refname}{\hfil ²Î \quad ¿¼\quad ÎÄ\quad Ï×}

%\end{CJK*}
\end{document}